\documentclass[draft,12pt, reqno]{amsart}
\usepackage{amssymb, amsmath}
\usepackage{url}
\usepackage{enumerate}
\usepackage{mathabx}

\newtheorem{theorem}{Theorem}[section]
\newtheorem{corollary}[theorem]{Corollary}

\theoremstyle{remark}

\theoremstyle{definition}

\newtheorem{example}[theorem]{Example}

\numberwithin{equation}{section}
\numberwithin{theorem}{section}

\newcommand{\N}{{\mathbb N}}
\newcommand{\R}{{\mathbb R}}

\newcommand{\fn}{\!:\!}

\newcommand{\fhat}{\hat{f}}

\providecommand{\abs}[1]{\lvert#1\rvert}
\providecommand{\norm}[1]{\lVert#1\rVert}

\providecommand{\conv}[2]{#1\ast#2}
\newcommand{\Lone}{L^1(\R)}
\newcommand{\intab}{\int_\alpha^\beta}

\newcommand{\intinf}{\int_{-\infty}^{\infty}}

\begin{document}
\date{Preprint February 1, 2022.  To appear in {\it Journal of Classical Analysis}.}
\title[Fourier transform inversion]
{Fourier transform inversion in the Alexiewicz norm}
\author{Erik Talvila}
\address{Department of Mathematics \& Statistics\\
University of the Fraser Valley\\
Abbotsford, BC Canada V2S 7M8}
\email{Erik.Talvila@ufv.ca}

\begin{abstract}
If $f\in L^1({\mathbb R})$ it is proved that $\lim_{S\to\infty}\lVert f-f\ast D_S\rVert=0$,
where $D_S(x)=\sin(Sx)/(\pi x)$ is the Dirichlet kernel and $\lVert f\rVert =
\sup_{\alpha<\beta}|\int_{\alpha}^{\beta}f(x)\,dx|$ is the Alexiewicz norm.  This gives a
symmetric inversion of the Fourier transform on the real line.  An asymmetric inversion
is also proved.  The results also hold for a measure given by $dF$ where $F$ is a continuous
function of bounded variation.  Such measures need not be absolutely continuous with respect to Lebesgue measure.
An example shows there is $f\in L^1({\mathbb R})$ such that $\lim_{S\to\infty}
\rVert f-f\ast D_S\lVert_1\neq 0$.
\end{abstract}

\maketitle

\section{Introduction}\label{sectionintroduction}
If $f\fn\R\to\R$ then its Fourier transform is $\fhat(s)=\intinf e^{-ist}f(t)\,dt$.  
Under the condition 
$f\in \Lone$ it is known that $\fhat$ is uniformly continuous on
$\R$.  The Riemann-Lebesgue lemma says that $\fhat$ vanishes at infinity.  
The inversion formula is $f(x)=\frac{1}{2\pi}\intinf e^{isx}\fhat(s)\,ds$.
Further assumptions, such as $\fhat\in\Lone$, are needed for this formula to hold.
See \cite{benedetto}  and \cite{titchmarsh} for such background information.

In this paper we prove the Fourier inversion theorem in the Alexiewicz norm, $\lim_{S\to\infty}\norm{
f-\conv{f}{D_S}}=0$.  The Alexiewicz norm is $\norm{f}=\sup_{\alpha<\beta}\abs{\int_\alpha^\beta f(x)\,dx}$.
It is useful for conditionally convergent integrals \cite{alexiewicz}, \cite{swartz}.  Note that $\norm{f}\leq\norm{f}_1$ with equality
if and only if $f$ is of one sign almost everywhere.  For $S>0$, the family of functions, $D_S(x)=\sin(Sx)/(\pi x)$,
is known as the Dirichlet kernel.  Notice that $\intinf D_S(x)\,dx=1$.

The Dirichlet kernel arises from a symmetric inversion integral.  Let $S>0$.
By the Fubini--Tonelli theorem,
\begin{eqnarray*}
\frac{1}{2\pi}\int_{-S}^S e^{isx}\fhat(s)\,ds & = & \frac{1}{2\pi}\intinf f(t)\int_{-S}^S e^{-i(t-x)s}\,ds\,dt\\
 & = & \frac{1}{\pi}\intinf f(t)\frac{\sin[(t-x)S]}{t-x}\,dt\\
 & = & f\ast D_S(x).
\end{eqnarray*}
The symmetric inversion formula is then  written as convolution with the Dirichlet kernel,
$$
f(x)=
\lim_{S\to\infty}\frac{1}{2\pi}\int_{-S}^S e^{isx}\fhat(s)\,ds =\lim_{S\to\infty} f\ast D_S(x).
$$
For this formula to hold pointwise, further conditions must be imposed on $f$ or $\fhat$.  For example,
it holds at points of differentiability of $f$ \cite{koekoek} or everywhere if $f$ is absolutely
continuous and $f'\in\Lone$ \cite{talvilafouriermaa}.

We prove in Theorem~\ref{theoreminversion} that the inversion holds in the Alexiewicz norm, $\lim_{S\to\infty}
\norm{f-f\ast D_S}=0$ for all $f\in\Lone$.  The proof is elementary and does not use any machinery from 
Fourier analysis.  For example, it does not use the Riemann--Lebesgue lemma.  
It should be noted that the same result is false for the $L^1$ norm,
i.e., there is $f\in\Lone$ such that $\lim_{S\to\infty}
\norm{f-f\ast D_S}_1\neq 0$.  See Example~\ref{exampleL1}, where $f\in\Lone$ is given so that
$\norm{f\ast D_S}_1$  exists for no $S>0$.

In Section~\ref{sectionalexnorminversion} we prove that if $F$ is a
continuous function of bounded variation and $f=dF$ is the associated signed Lebesgue--Stieltjes measure then
the inversion formula holds for $dF$.  Such measures need not be absolutely continuous with respect to
Lebesgue measure.  An example shows that inversion can fail when $F$ is not continuous.

In Section~\ref{sectionasymmetric}  we prove the inversion still holds with an 
asymmetric
inversion integral.

A similar result holds for Fourier series \cite{talvilafourierseries}.

\section{Alexiewicz norm inversion theorem}\label{sectionalexnorminversion}
\begin{theorem}\label{theoreminversion}
Let $f\in\Lone$.   Then $\lim_{S\to\infty}
\norm{f-f\ast D_S}=0$, where $D_S$ is the Dirichlet kernel.
\end{theorem}
\begin{proof}
Let $-\infty<\alpha<\beta<\infty$ and let $S>0$. Let $F(x)=\int_{-\infty}^xf(t)\,dt$.  By the Fubini--Tonelli theorem,
\begin{align*}
&\intab[f(x)-f\ast D_S(x)]\,dx  =  \intab\intinf\left[f(x) D_1(t)-
f(x-t)D_S(t)\right]\,dt\,dx\\
&=  \frac{1}{\pi}\intinf\frac{\sin(t)}{t}\intab[f(x)-f(x-t/S)]\,dx\,dt\\
&=\frac{1}{\pi}\intinf\frac{\sin(t)}{t}\left[F(\beta)-F(\beta-t/S)-F(\alpha)+F(\alpha-t/S)\right]\,dt.
\end{align*}
Let $T>0$.
We split this last integral into integration over the intervals $(-T,T)$, $(T,\infty)$ and
$(-\infty,-T)$.  The supremum over $\alpha<\beta$ of the absolute value of each such 
integral is then shown to have limit $0$ as $S\to\infty$.

First note that
\begin{align*}
&\left|\int_{-T}^T\frac{\sin(t)}{t}\left[F(\beta)-F(\beta-t/S)-F(\alpha)+F(\alpha-t/S)\right]\,dt\right|\\
&\leq 2T\left(\sup_{\substack{\beta\in\R\\\abs{t}\leq T}}\abs{F(\beta)-F(\beta-t/S)}
+\sup_{\substack{\alpha\in\R\\\abs{t}\leq T}}\abs{F(\alpha)-F(\alpha-t/S)}\right)\\
&\to 0 \text{ as } S\to\infty \text{ since } F \text{ is uniformly continuous on } \R.
\end{align*}

Let $\epsilon>0$.
Let $g(t)=\sin(t)/t$.
Note that $\lim_{t\to\infty}F(\beta-t/S)=0$.
Note that $\int_T^\infty\abs{f(\beta-t/S)}\,dt/S=\int_{-\infty}^{\beta-T/S}\abs
{f(u)}\,du\leq \norm{f}_1$ and $\lim_{T\to\infty}\int_T^\infty g(t)\,dt=0$.
Now integrate by parts to get
\begin{align*}
&\left|\int_T^\infty g(t)\left[F(\beta)-F(\beta-t/S)-F(\alpha)+F(\alpha-t/S)\right]\,dt\right|\\
&=\left|[F(\beta)-F(\alpha)]\int_T^\infty g(t)\,dt-\int_T^\infty\int_T^tg(u)\,du[f(\beta-t/S)
-f(\alpha-t/S)]\,dt/S\right|\\
&\leq 3\norm{f}_1\norm{\chi_{(T,\infty)}g}<\epsilon \text{ for large enough } T.
\end{align*}
Similarly for integration over the interval $(-\infty,-T)$.
\end{proof}

The essential parts of the proof of the theorem use the Fubini--Tonelli theorem and the fact that $F$
is uniformly continuous and
has finite variation.
Thus, the theorem extends to measures that arise from continuous functions of bounded variation.
If $\mu$ is a signed measure then its Alexiewicz norm is
$$
\norm{\mu}=\sup_{\alpha<\beta}\left|\int_{(\alpha,\beta)}d\mu\right|=\sup_{\alpha<\beta}
\abs{\mu((\alpha,\beta))}.
$$

\begin{corollary}\label{corollarymeasure}
Let $F\in C(\R)$ such that $F$ is of bounded variation.  Define a signed measure $\mu_F=dF$.
Then $\lim_{S\to\infty}\norm{\mu_F-\mu_F\ast D_S}=0$.
\end{corollary}

If $F$ is of bounded variation then $\lim_{x\to\infty}F(x)$ and $\lim_{x\to-\infty}F(x)$
exist so that if $F(\infty)$ and $F(-\infty)$ are defined with these respective limits then
$F$ is continuous on the extended real line $[-\infty,\infty]$.

The continuity condition in Corollary~\ref{corollarymeasure} cannot be dropped.
\begin{example}
The Dirac measure, $\delta$, is generated by the Heaviside step function, $H(x)=0$ for $x<0$ and
$H(x)=1$ for $x>0$.  We have $\hat{\delta}=1$.  This gives 
$$
\frac{1}{2\pi}\int_{-S}^S e^{isx}\hat{\delta}(s)\,ds  =  \frac{1}{2\pi}\int_{-S}^S e^{isx}\,ds
  =  D_S(x).
$$

Let $0<\alpha<\beta$.  Then
$$
\left|\intab d\delta(x)-\frac{1}{\pi}\intab\frac{\sin(Sx)}{x}\,dx\right|=\frac{1}{\pi}\left|\intab
\frac{\sin(Sx)}{x}\,dx\right|=\frac{1}{\pi}\left|\int_{\alpha S}^{\beta S}
\frac{\sin(x)}{x}\,dx\right|.
$$
If we let $\alpha=1/S^2$ and $\beta=1$ then we see that $\limsup_{S\to\infty}\norm{
\delta-\delta\ast D_S}\geq 1/2$.
\end{example}

\begin{example}\label{exampleL1}
Let $a>0$ and take $f=\chi_{(0,a)}$.  For $x>0$, integrate by parts to get
$$
\pi f\ast D_S(x)  =  \int_{(x-a)S}^{xS}\sin(t)\,\frac{dt}{t}
  =  \frac{\cos[(x-a)S]}{(x-a)S}-\frac{\cos(xS)}{xS}-\int_{(x-a)S}^{xS}\cos(t)\,
\frac{dt}{t^2}.
$$
Note that $1/(x-a)=1/x+a/x^2+O(1/x^3)$ as $x\to\infty$.  Then
$$
\left|\int_{(x-a)S}^{xS}\cos(t)\,\frac{dt}{t^2}\right|\leq 
\left|\int_{(x-a)S}^{xS}\frac{dt}{t^2}\right|
=\frac{1}{(x-a)S}-\frac{1}{xS}=\frac{a}{x^2S}+O(1/x^3) \text{ as } x\to\infty.
$$
This gives
$$
\pi f\ast D_S(x)
  =  \frac{\sin(aS)\sin(xS)-[1-\cos(aS)]\cos(xS)}{xS} +O\left(\frac{1}{x^2}\right)
\text{ as } x\to\infty.
$$
Hence, if $S\neq 2n\pi/a$ for some $n\in\N$ then $f\ast D_S\not\in \Lone$.

If we write $f_a=\chi_{(0,a)}$ then 
\begin{align*}
&\pi(f_a+f_b)\ast D_S(x)\\
&=  \frac{[\sin(aS)+\sin(bS)]\sin(xS)-[2-\cos(aS)-\cos(bS)]\cos(xS)}{xS} 
+O\left(\frac{1}{x^2}\right).
\end{align*} 
It follows that $(f_a+f_b)\ast D_S\in L^1(\R)$ if and only if $\sin(aS)=\sin(bS)=0$
and $\cos(aS)=\cos(bS)=1$.  Taking $a=2\pi$ and $b=2$ then gives an example of a 
function $f\in\Lone$ such that $\norm{f-f\ast D_S}_1$ fails
to exist for each $S>0$.
\end{example}

\section{Asymmetric kernel}\label{sectionasymmetric}
Let $S_1,S_2>0$.  Consider the asymmetric inversion,
\begin{eqnarray*}
\frac{1}{2\pi}\int_{-S_1}^{S_2}e^{isx}\fhat(s)\,ds & = & \frac{1}{2\pi}\intinf f(t)
\int_{-S_1}^{S_2}e^{-i(t-x)s}\,ds\,dt\\
 & = & \frac{i}{2\pi}\intinf f(t)\left[e^{-i(t-x)S_2}-e^{i(t-x)S_1}\right]\frac{dt}{t-x}\\
 & = & f\ast A_{S_1,S_2}(x),
\end{eqnarray*}
where the asymmetric kernel is 
$$
A_{S_1,S_2}(x)   =  \frac{e^{ixS_2}-e^{-ixS_1}}{2\pi ix}
  =   \frac{D_{S_1}(x)}{2}+
\frac{D_{S_2}(x)}{2}+ B_{S_1,S_2}(x),
$$
where
$$
B_{S_1,S_2}(x)=\frac{1}{\pi ix}\sin\left[\left(\frac{S_1+S_2}{2}\right)x\right]
\sin\left[\left(\frac{S_1-S_2}{2}\right)x\right].
$$
Notice that $\intinf A_{S_1,S_2}(x)\,dx=1$.  As with Theorem~\ref{theoreminversion} there
is inversion in the Alexiewicz norm, provided the ratios $S_1/S_2$ and $S_2/S_1$ remain bounded.

\begin{theorem}\label{theoremasymmetric}
Let $q\geq p>0$.  Let $f\in\Lone$. Let $S_1,S_2\to\infty$ such that $p\leq S_1/S_2\leq q$.
Then $\norm{f-f\ast A_{S_1,S_2}}\to 0$.
\end{theorem}
\begin{proof}
The proof is similar to that of Theorem~\ref{theoreminversion}.
Let $-\infty<\alpha<\beta<\infty$ and let $S_1,S_2>0$. Let $F(x)=\int_{-\infty}^xf(t)\,dt$.
Write
$$
f-f\ast A_{S_1,S_2}=\frac{1}{2}\left(f-f\ast D_{S_1}\right)
+\frac{1}{2}\left(f-f\ast D_{S_1}\right)-f\ast B_{S_1,S_2}.
$$
Due to Theorem~\ref{theoreminversion} we just need to consider the last term above.
Since $B_{S_1,S_2}$ is an odd function we can write
\begin{align*}
&-\pi i\intab f\ast B_{S_1,S_2}(x)\,dx\\
&=\pi i\intab f(x)\,dx\intinf B_{S_1,S_2}(t)\,dt
-\pi i\intab\intinf f(x-t) B_{S_1,S_2}(t)\,dt\,dx\\
&=\intab\intinf\left[f(x)-f(x-2t/(S_1+S_2))\right]
\sin(t)\sin\left[\left(\frac{S_1-S_2}{S_1+S_2}\right)t\right]\frac{dt}{t}\,dx,
\end{align*}
after the change of variables $t\mapsto 2t/(S_1+S_2)$.

Let $T>0$.  Then
\begin{align*}
&\left|\intab\int_{-T}^T\left[f(x)-f(x-2t/(S_1+S_2))\right]
\sin(t)\sin\left[\left(\frac{S_1-S_2}{S_1+S_2}\right)t\right]\frac{dt}{t}\,dx\right|\\
&\leq 2T\left(\sup_{\substack{\beta\in\R\\\abs{t}\leq T}}\abs{F(\beta)-F(\beta-2t/(S_1+S_2))}
+\sup_{\substack{\alpha\in\R\\\abs{t}\leq T}}\abs{F(\alpha)-F(\alpha-2t/(S_1+S_2))}\right)\\
&\to 0 \text{ as } S_1,S_2\to\infty \text{ since } F \text{ is uniformly continuous on } \R.
\end{align*}

And,
\begin{align*}
&\intab\int_T^\infty \left[f(x)-f(x-2t/(S_1+S_2))\right]\sin(t)\sin\left[\left(\frac{S_1-S_2}{S_1+S_2}\right)t\right]\frac{dt}{t}
\,dx\\
&=\frac{1}{2}\int_T^\infty \left[F(\beta)-F(\beta-2t/(S_1+S_2))\right]
\left[\cos\left(\frac{2S_2t}{S_1+S_2}\right)-\cos\left(\frac{2S_1t}{S_1+S_2}\right)\right]\,\frac{dt}{t}\\
&\quad + \text{ similar terms involving } \alpha.
\end{align*}
Let $\epsilon>0$ be given.
Let $h(t)=\cos(t)/t$.  Write $A_1=2S_1T/(S_1+S_2)$ and $A_2=2S_2T/(S_1+S_2)$.  Use
the change of variables $u_2=2S_2t/(S_1+S_2)$ and $u_1=2S_1t/(S_1+S_2)$.
Integrate by parts as in the last paragraph of the proof of Theorem~\ref{theoreminversion}.
Then
\begin{align*}
&\int_T^\infty \left[F(\beta)-F(\beta-2t/(S_1+S_2))\right]
\left[\cos\left(\frac{2S_2t}{S_1+S_2}\right)-\cos\left(\frac{2S_1t}{S_1+S_2}\right)\right]\,\frac{dt}{t}\\
&=\int_{A_2}^\infty \left[F(\beta)-F(\beta-u_2/S_2)\right]h(u_2)\,du_2
-\int_{A_1}^\infty \left[F(\beta)-F(\beta-u_1/S_1)\right]h(u_1)\,du_1\\
&=F(\beta)\int_{A_2}^\infty h(u_2)\,du_2-\int_{A_2}^\infty \int_{A_2}^{u_2}h(v)\,dv 
\,f(\beta-u_2/S_2)\frac{du_2}{S_2}\\
&\quad-F(\beta)\int_{A_1}^\infty h(u_1)\,du_1+\int_{A_1}^\infty \int_{A_1}^{u_1}h(v)\,dv 
\,f(\beta-u_1/S_1)\frac{du_1}{S_1}.
\end{align*}
Notice that $A_2=2S_2T/(S_1+S_2)\geq 2T/(q+1)$ and $S_1=2S_1T/(S_1+S_2)\geq 2T/(1+1/p)$.
Hence, if $T$ is large then so are $A_2$ and $A_1$.
This gives
\begin{align*}
&\left|\int_T^\infty \left[F(\beta)-F(\beta-2t/(S_1+S_2))\right]
\left[\cos\left(\frac{2S_2t}{S_1+S_2}\right)-\cos\left(\frac{2S_1t}{S_1+S_2}\right)\right]\,\frac{dt}{t}\right|\\
&\leq 2\norm{f}_1\left(\norm{\chi_{(A_2,\infty)}h}+\norm{\chi_{(A_1,\infty)}h}\right)\\
&<\epsilon \text{ for large enough } T.
\end{align*}
The same estimates hold for the terms containing $\alpha$.
\end{proof}

\end{document}